\documentclass[11pt]{amsart}
\usepackage{amssymb, latexsym}
\theoremstyle{plain}
\newtheorem{theorem}{Theorem}
\newtheorem{corollary}{Corollary}
\newtheorem*{thm-cheb}{Theorem (Chebyshev)}

\newtheorem{proposition}{Proposition}

\newtheorem*{2'}{Theorem 2'}
\newtheorem*{3'}{Theorem 3'}

\theoremstyle{remark}

\newtheorem*{Remark 1}{Remark 1}
\newtheorem*{Remark 2}{Remark 2}
\newtheorem*{Remark 3}{Remark 3}
\newtheorem*{Remark 4}{Remark 4}

\numberwithin{equation}{section}

\begin{document}

\title[  Clustering under  Mallows  and  general $p$-shifted distributions]
 {Clustering of consecutive numbers in permutations under Mallows distributions  and super-clustering under general $p$-shifted distributions}

\author{Ross G. Pinsky}


\address{Department of Mathematics\\
Technion---Israel Institute of Technology\\
Haifa, 32000\\ Israel}
\email{ pinsky@math.technion.ac.il}

\urladdr{http://www.math.technion.ac.il/~pinsky/}

\date{}


\begin{abstract}
Let $A^{(n)}_{l;k}\subset S_n$ denote the set of permutations of $[n]$ for which the set of $l$ consecutive numbers
$\{k, k+1,\cdots, k+l-1\}$ appears
in a  set  of consecutive positions.
Under the uniformly probability measure $P_n$ on $S_n$,
one has $P_n(A^{(n)}_{l;k})\sim\frac{l!}{n^{l-1}}$ as $n\to\infty$. In one part of this paper we consider the probability of clustering of consecutive numbers
under Mallows distributions $P_n^q$, $q>0$. Because of a duality, it suffices to consider $q\in(0,1)$.
We show that for $q_n=1-\frac c{n^\alpha}$, with $c>0$ and
 $\alpha\in(0,1)$,   $P_n^q(A^{(n)}_{l;k_n})$ is on  the order $\frac1{n^{\alpha(l-1)}}$, uniformly over all sequences $\{k_n\}_{n=1}^\infty$. Thus, letting
$N^{(n)}_l=\sum_{k=1}^{n-l+1}1_{A^{(n)}_{l;k}}$ denote the number of sets of $l$ consecutive numbers
appearing in sets of consecutive positions, we have
 \begin{equation*}
\lim_{n\to\infty}E_n^{q_n}N^{(n)}_l=\begin{cases}\infty,\ \text{if}\ l<\frac{1+\alpha}\alpha;\\ 0,\ \text{if} \ l>\frac{1+\alpha}\alpha.
\end{cases}.
\end{equation*}
We also consider the cases $\alpha=1$ and $\alpha>1$.
In the other part of the paper we consider general $p$-shifted distributions, of which the Mallows distribution is a particular case. We calculate explicitly the quantity
$\lim_{l\to\infty}\liminf_{n\to\infty}P_n^q(A^{(n)}_{l;k_n})=\lim_{l\to\infty}\limsup_{n\to\infty}P_n^q(A^{(n)}_{l;k_n})$ in terms of the $p$-distribution. When this quantity is positive, we say that super-clustering occurs.
In particular, super-clustering occurs for the Mallows distribution with parameter $q\neq1$.
We also give a new characterization of $p$-shifted distributions.
\end{abstract}

\maketitle
\section{Introduction and Statement of Results}
Let $l\ge 2$ be an integer.
Let $P_n$ denote the uniform probability measure on the set $S_n$ of permutations of $[n]:=\{1,\cdots, n\}$, and  denote a permutation $\sigma\in S_n$ by $\sigma=\sigma_1\sigma_2\cdots \sigma_n$.
The set of  $l$ consecutive numbers $\{k,k+1,\cdots, k+l-1\}\subset [n]$ appears
in a  set  of consecutive positions in the permutation if there exists an $m$ such that $\{k,k+1,\cdots, k+l-1\}=\{\sigma_m,\sigma_{m+1},\cdots, \sigma_{m+l-1}\}$.
Let $A^{(n)}_{l;k}\subset S_n$ denote the event that
 the set of  $l$ consecutive numbers $\{k, k+1,\cdots, k+l-1\}$ appears
in a  set  of consecutive positions.
It is immediate that for any $1\le k,m\le n-l+1$,
the probability that  $\{k,k+1,\cdots, k+l-1\}=\{\sigma_m,\sigma_{m+1},\cdots, \sigma_{m+l-1}\}$ is equal to $\frac{l!(n-l)!}{n!}$.
Thus,
\begin{equation}\label{uniformprob}
P_n(A^{(n)}_{l;k})=(n-l+1)\frac{l!(n-l)!}{n!}\sim\frac{l!}{n^{l-1}},\ \text{as}\ n\to\infty, \text{for}\ l\ge2.
\end{equation}

Let $A^{(n)}_l=\cup_{k=1}^{n-l+1}A^{(n)}_{l;k}$ denote the event that there exists a set of   $l$ consecutive numbers
appearing in a  set of consecutive positions, and let
$N^{(n)}_l=\sum_{k=1}^{n-l+1}1_{A^{(n)}_{l;k}}$ denote the number of sets of $l$ consecutive numbers
appearing in sets of consecutive positions.
Then
\begin{equation}\label{uniformexpect}
E_nN^{(n)}_l=(n-l+1)^2\frac{l!(n-l)!}{n!}\sim\frac{l!}{n^{l-2}},\ \text{as}\ n\to\infty,\ \text{for}\ l\ge2.
\end{equation}
Using the inequality
$$
\sum_{k=1}^{n-k+1}P_n(A^{(n)}_{l;k})-\sum_{1\le j<k\le n-l+1}P_n(A^{(n)}_{l;j}\cap A^{(n)}_{l;k})\le P_n(A^{(n)}_l)\le\sum_{k=1}^{n-k+1}P_n(A^{(n)}_{l;k}),
$$
along with the fact that for $j,k,m,r$,
 with
$\{j,j+1,\cdots, j+l-1\}\cap\{k,k+1,\cdots, k+l-1\}=\emptyset$ and $\{m,m+1,\cdots, m+l-1\}\cap\{r,r+1,\cdots, r+l-1\}=\emptyset$,
the probability that both $\{k,k+1,\cdots, k+l-1\}=\{\sigma_m,\sigma_{m+1},\cdots, \sigma_{m+l-1}\}$
and $\{j,j+1,\cdots, j+l-1\}=\{\sigma_r,\sigma_{r+1},\cdots, \sigma_{r+l-1}\}$ is equal to,
$\frac{(l!)^2(n-2l)!}{n!}$, it is easy to show that
\begin{equation}\label{uniformtotalprob}
P_n(A^{(n)}_l)\sim\frac{l!}{n^{l-2}},\ \text{as}\ n\to\infty,\ \text{for}\ l\ge3.
\end{equation}

It follows from \eqref{uniformexpect}  (or from \eqref{uniformtotalprob}) that for $l\ge3$, the sequence $\{N_l^{(n)}\}_{n=1}^\infty$ converges to zero in probability. On the other hand,
 when $l=2$, $\{N_l^{(n)}\}_{n=1}^\infty$ converges in distribution to a Poisson random variable with parameter $2$. This result goes back over 75 years; see
\cite{W}, \cite{K}.

In one of the two  parts of this paper, we obtain results in the spirit of \eqref{uniformprob} and \eqref{uniformexpect}       in the case that the uniform probability measure $P_n$ is replaced by
the Mallows measure $P_n^{q_n}$, with $q_n\to1$ at various rates.
The Mallows measures $P_n^q$ are described below.
The Mallows measure with $q=1$ is the uniform measure.

For fixed $q\neq1$, it turns out
that $P_n^q(A^{(n)}_{l;k})$ remains bounded away from 0 as $n\to\infty$, for all $l$. In the other part of this paper we consider so-called $p$-shifted distributions $P_n^{(\{p_j\}_{j=1}^\infty)}$ on $S_n$, of which the Mallows measure $P_n^q$ is a particular example. Here $\{p_j\}_{j=1}^\infty$, with $p_j>0$, for all $j$, is a probability distribution on $\mathbb{N}$: $\sum_{j=1}^\infty p_j=1$.
We calculate $\lim_{l\to\infty}\lim_{n\to\infty}P_n^{(\{p_j\}_{j=1}^\infty)}(A^{(n)}_{l;k})$ explicitly. This reveals a necessary and sufficient condition on the distribution $\{p_j\}_{j=1}^\infty$ in order for the above limit to be positive. When this limit is positive, we  say that \it super-clustering\rm\ occurs.
In particular,  super-clustering occurs for the $q$-Mallows distributions for $q\neq1$.
We also give a new characterization of $p$-shifted measures, which may be of some independent interest.

We turn now to a description of the results of the   part of the paper concerning specifically the Mallows distributions.

\medskip

\bf\noindent The behavior of the probability of $A^{(n)}_{l;k}$
under Mallows distributions.\rm\
 For each $q>0$, the Mallows distribution with parameter $q$ is the probability measure $P_n^q$ on $S_n$ defined by
\begin{equation}\label{mallowsdef}
P_n^q(\sigma)=\frac{q^{\text{inv}(\sigma)}}{Z_n(q)}, \sigma\in S_n,
\end{equation}
where $\text{inv}(\sigma)$ is the number of inversions in $\sigma$, and $Z_n(q)$ is the normalization constant, given by
$$
Z_n(q)=\prod_{k=1}^n\frac{1-q^k}{1-q}.
$$
Thus, for $q\in(0,1)$, the distribution favors permutations with few inversions, while for $q>1$, the distribution favors permutations with many inversions. Of course, the case $q=1$ yields the uniform distribution.
Recall that the \it reverse\rm\  of a permutation $\sigma=\sigma_1\cdots\sigma_n$ is the permutation $\sigma^{\text{rev}}:=\sigma_n\cdots\sigma_1$.
The Mallows distributions satisfy the following duality between $q>1$ and $q\in(0,1)$:
$$
P_n^q(\sigma)=P_n^{\frac1q}(\sigma^{\text{rev}}),\ \text{for}\ q>0, \sigma\in S_n\ \text{and}\ n=1,2,\cdots.
$$
Since the set $A^{(n)}_{l;k}$   is invariant under reversal, for our study of clustering it suffices to consider the case that $q\in(0,1)$.

When $q\to0$, the Mallows distribution $P_n^q$ converges weakly to the degenerate distribution on
the identity permutation, and of course the identity permutation belongs to $A^{(n)}_{l;k}$ for all $k$ and $l$. Because the smaller $q$ is, the more the distribution favors permutations with few inversions, and as such, the smaller $q$ is, the more the
 distribution favors permutations which are close to the identity permutation,
 it seems intuitive that the smaller $q$ is, the more clustering there will be.
 However, whereas the structure of the Mallows distribution lends itself
 naturally  to  proving  theorems concerning the inversion statistic
\cite{P}, it is less transparent how to exploit that structure   with regard to this clustering statistic.
For example, the set $A^{(n)}_{l;k}$ is the disjoint union of the $n-l+1$ sets
$\{k,k+1,\cdots, k+l-1\}=\{\sigma_m,\sigma_{m+1},\cdots, \sigma_{m+l-1}\},\ m=1,\cdots, n-l+1$.
In the case of the uniform distribution, these $n-l+1$ sets all have the same probability. However, in the case of
$P_n^q$, $q\in(0,1)$, we expect that for certain $m$, these sets will have probability less than what they have under the uniform
distribution, and for other $m$ these sets
will have probability greater than what they have under the uniform distribution.

For results concerning the behavior under a Mallows distribution of other permutation statistics, such as cycle counts and increasing subsequences, see \cite{BB}, \cite{BP} and \cite{GP}.

Our first theorem gives   asymptotic results in the case that $q=q_n=1-\frac c{n^\alpha}$ with $c>0$ and $\alpha\in(0,1)$.
We use the notation $a_n\lesssim b_n$ as $n\to\infty$ to indicate that $\limsup_{n\to\infty}\frac{a_n}{b_n}\le 1$.
\pagebreak
\begin{theorem}\label{maintheorem}
Let $A^{(n)}_{l;k}\subset S_n$ denote the event that  the set of $l$ consecutive numbers $\{k,k+1,\cdots, k+l-1\}$
 appears in a set of $l$ consecutive positions.
Let $q_n=1-\frac c{n^\alpha}$, with $c>0$ and $\alpha\in(0,1)$.
Then
\begin{equation}\label{mainresultgen}
\frac{\big((l-1)!\big)^2}{(2l)!}\thinspace \frac{c^{\thinspace l-1}\thinspace l!} {n^{\alpha(l-1)}}
\lesssim P_n^{q_n}(A^{(n)}_{l;k_n})\lesssim \frac1l\thinspace \frac{c^{\thinspace l-1}\thinspace l!} {n^{\alpha(l-1)}},
\end{equation}
for any choice of $\{k_n\}_{n=1}^\infty$, and the asymptotics  are uniform over all $\{k_n\}_{n=1}^\infty$.
If $k_n$ satisfies $\frac{\min(k_n,n-k_n)}{n^\alpha}=\infty$,
then \eqref{mainresultgen} holds with an improved upper bound:
\begin{equation}\label{mainresult}
\frac{\big((l-1)!\big)^2}{(2l)!}\thinspace \frac{c^{\thinspace l-1}\thinspace l!} {n^{\alpha(l-1)}}
\lesssim P_n^{q_n}(A^{(n)}_{l;k_n})\lesssim \big(\int_0^1x^{l-1}e^{-(l-1)x}dx\big)\thinspace \frac{c^{\thinspace l-1}\thinspace l!} {n^{\alpha(l-1)}}.
\end{equation}
\end{theorem}
\medskip
Recall that $N^{(n)}_l=\sum_{k=1}^{n-l+1}1_{A^{(n)}_{l;k}}$ denotes the number of sets of $l$ consecutive numbers
appearing in sets of consecutive positions.
Theorem \ref{maintheorem} yields the following corollary.
\begin{corollary}
Let $q_n=1-\frac c{n^\alpha}$ with $c>0$ and $\alpha\in(0,1)$,
Then there exist constants $C_l^{(-)},C_l^{(+)}>0$ such that
\begin{equation*}
C_l^{(-)} n^{1-(l-1)\alpha}\le E_n^{q_n}N^{(n)}_l\le C_l^{(+)} n^{1-(l-1)\alpha}.
\end{equation*}
In particular,
\begin{equation*}
\lim_{n\to\infty}E_n^{q_n}N^{(n)}_l=\begin{cases}\infty,\ \text{if}\ l<\frac{1+\alpha}\alpha;\\ 0,\ \text{if} \ l>\frac{1+\alpha}\alpha.
\end{cases}.
\end{equation*}
\end{corollary}

\bf\noindent Remark 1.\rm\ For $\tau\in S_l$, let
$A^{(n)}_{l,\tau;k}\subset A^{(n)}_{l;k}$ denote the event that
the set of  $l$ consecutive numbers $\{k,k+1,\cdots, k+l-1\}\subset [n]$ appears
in a  set  of consecutive positions in the permutation and also that the relative positions of these consecutive numbers  correspond to the
permutation $\tau$. That is,  $\{k,k+1,\cdots, k+l-1\}=\{\sigma_m,\sigma_{m+1},\cdots, \sigma_{m+l-1}\}$,
for some $m$, and $\sigma_{m+i-1}-(k-1)=\tau_i,\ i=1,\cdots, l$. Then $A^{(n)}_{l;k}=\cup_{\tau\in S_l}A^{(n)}_{l,\tau;k}$.
Small changes in the proof of Theorem \ref{maintheorem}, which   we leave to the reader, show that
\eqref{mainresultgen} and \eqref{mainresult} hold with
$P_n^{q_n}(A^{(n)}_{l;k_n})$
replaced by $P_n^{q_n}(A^{(n)}_{l,\tau;k_n})$ and with $l!$ deleted from the numerator in the upper and lower bounds, for all $\tau\in S_l$.
In particular, if $\tau=id$, then $A^{(n)}_{l,\tau;k_n}$ is the event that the numbers $\{k,\cdots, k+l-1\}$
form an \it increasing run\rm\ in the permutation, and if $\tau$ satisfies $\tau^{\text{rev}}=id$, then
$A^{(n)}_{l,\tau;k_n}$ is the event that the numbers $\{k,\cdots, k+l-1\}$
form a \it decreasing run\rm\ in the permutation.
\medskip

\bf\noindent Remark 2.\rm\
Let $K^{(-)}(l)=\frac{\big((l-1)!\big)^2}{(2l)!}$ and $K^{(+)}(l)=\int_0^1x^{l-1}e^{-(l-1)x}dx$ denote the coefficients
of $\frac{\thinspace c^{\thinspace l-1}\thinspace l!} {n^{\alpha(l-1)}}$
on the left and right hand sides respectively of \eqref{mainresult}.
We have $K^{(-)}(l)\sim\sqrt{\pi}\ l^{-\frac32}4^{-l}$ as $l\to\infty$.
One can show that
$$
K^{(+)}(l)=\int_0^1x^{l-1}e^{-(l-1)x}dx=\frac{(l-1)!}{(l-1)^l}\big(1-e^{-(l-1)}\sum_{i=0}^{l-1}\frac{(l-1)^i}{i!}\big).
$$
Thus, $K^{(+)}(l)\lesssim \frac{(l-1)!}{(l-1)^l}\sim \sqrt{2\pi}\thinspace e\ l^{-\frac12}e^{-l}$, as $l\to\infty$.
On the other hand, a rudimentary asymptotic analysis we performed on  the interval $[\frac{l-1}l-l^{-\frac12},1]$ yields
$K^{(+)}(l)\gtrsim e^{\frac12}\ l^{-\frac12}e^{-l}$, as $l\to\infty$.
We have $K^{(-)}(2)=\frac1{12}\approx0.083$ and $K^{(+)}(2)=1-\frac2e\approx0.281$.
\medskip

Now we consider the cases $q=q_n=1-\frac cn$ and $q=q_n=1-o(\frac1n)$.
\begin{theorem}\label{alpha1}
Let $A^{(n)}_{l;k}\subset S_n$ denote the event that  the set of $l$ consecutive numbers $\{k,k+1,\cdots, k+l-1\}$
 appears in a set of $l$ consecutive positions.

\noindent i.
Let $q_n=1-\frac cn$, with $c>0$. Let $k_n\sim dn$ with $d\in(0,1)$.
Then
\begin{equation}\label{alpha1-upper}
P_n^{q_n}(A^{(n)}_{l;k_n})\lesssim  \frac1{(1-e^{-cd})^l}\big(\int_{e^{-cd}}^1y^{l-1}
e^{(\log\frac{1-e^{-cd}}{1-e^{-c}})e^{cd}(l-1)y}dy\big)
 \thinspace \frac{c^{\thinspace l-1}\thinspace l!} {n^{(l-1)}}.
\end{equation}

\noindent ii. Let $q_n=1-o(\frac1n)<1$. Then for any choice of $\{k_n\}_{n=1}^\infty$,
\begin{equation}\label{littleo}
P_n^{q_n}(A^{(n)}_{l;k_n})\lesssim \frac{l!}{n^{l-1}}.
\end{equation}
\end{theorem}
\medskip

\bf\noindent Remark.\rm\ In part (i), we expect that the asymptotic behavior of $P_n^{q_n}(A^{(n)}_{l;k_n})$, when
$k_n\sim dn$, is in fact independent of $d\in(0,1)$.
We note for all $d\in(0,1)$,  the expression
$
\frac1{(1-e^{-cd})^l}\big(\int_{e^{-cd}}^1y^{l-1}
e^{(\log\frac{1-e^{-cd}}{1-e^{-c}})e^{cd}(l-1)cy}dy\big)c^{l-1}
$,
which multiplies $\frac{l!}{n^{l-1}}$ on the right hand side of \eqref{alpha1-upper},
converges to 1 when $c\to0$, thus matching up with \eqref{littleo}.
In the case of the uniform distribution ($q=1$), we have from \eqref{uniformprob} that
$P_n^1(A^{(n)}_{l;k_n})\sim\frac{l!}{n^{l-1}}$, for any choice of $k_n$. Since we expect
$P_n^q(A^{(n)}_{l;k_n})$ to be decreasing in $q$, we expect that
the asymptotic inequality in \eqref{littleo} is an asymptotic  equality.
That is, asymptotically, we expect that the cluster event $A^{(n)}_{l;k_n}$ cannot be used to distinguish between
$P_n^1$ and $P_n^{q_n}$, if $q_n=1-o(\frac1n)$.

\medskip

The following duality will be useful in the proofs of the above theorems.
Its short proof is given at the end of this section.
\begin{proposition}\label{dualprop}
\begin{equation}\label{duality}
P_n^q(A^{(n)}_{l;k})=P_n^q(A^{(n)}_{l;n+2-k-l}), \ k=1,2\dots, n-l+1.
\end{equation}
\end{proposition}

We now turn to a description of the results of the  other part of the paper, concerning $p$-shifted random permutations.
\medskip

\noindent \bf $p$-shifted distributions and super-clustering.\rm\ Denote by $S_\infty$ the set of permutations of $\mathbb{N}$. We build random permutations in $S_\infty$ and then project them down in a natural way to
S$_n$.
Let $p:=\{p_j\}_{j=1}^\infty$ be a probability distribution on $\mathbb{N}$ whose support is all of $\mathbb{N}$; that is,
$p_j>0$, for all $j\in\mathbb{N}$.  Take a countably infinite sequence of independent samples from this distribution: $n_1,n_2,\cdots$.
Now construct a random permutation $\Pi\in S_\infty$ as follows.
Let  $\Pi_1=n_1$ and then for $k\ge2$, let
  $\Pi_k=\psi_k(n_k)$, where $\psi_k$ is the increasing bijection from $\mathbb{N}$ to $\mathbb{N}-\{\Pi_1,\cdots,\Pi_{k-1}\}$.
Thus, for example, if the sequence of samples $\{n_j\}_{j=1}^\infty$ begins with  $7,3,4,3,7,2,1$, then the construction
yields the permutation $\Pi$ beginning with
$\Pi_1=7,\Pi_2=3,\Pi_3=5, \Pi_4=4,\Pi_5=11,\Pi_6=2,\Pi_7=1$.
The probability measure $P^{(\{p_j\}_{j=1}^\infty)}$ on $S_\infty$ is then the distribution of this random permutation $\Pi$.
We call $P^{(\{p_j\}_{j=1}^\infty)}$  the \it $p$-shifted distribution\rm\ and $\Pi$ a \it $p$-shifted random permutation\rm\ on $S_\infty$.

For $\sigma\in S_\infty$, write $\sigma=\sigma_1\sigma_2\cdots$.
For $n\in\mathbb{N}$,
define $\text{proj}_n(\sigma)\in S_n$ to be the permutation obtained from $\sigma$ by deleting $\sigma_i$ for
all $i$ satisfying $\sigma_i>n$. Thus, for $n=4$ and $\sigma=2539461\cdots$, one has $\text{proj}_4(\sigma)=2341$.
Given the $p$-shifted random permutation
$\Pi\in S_\infty$ that was constructed in   the previous paragraph, define
  $P_n^{(\{p_j\}_{j=1}^\infty)}$  as the distribution of the random permutation $\text{proj}_n(\Pi)$.
Equivalently,
given the probability measure $P^{(\{p_j\}_{j=1}^\infty)}$ on $S_\infty$ defined in the previous paragraph,
define the probability measure $P_n^{(\{p_j\}_{j=1}^\infty)}$ on $S_n$ by
$P_n^{(\{p_j\}_{j=1}^\infty)}(\sigma)=P^{(\{p_j\}_{j=1}^\infty)}(\text{proj}_n^{-1}(\sigma)),\ \sigma\in S_n$.
We call $P_n^{(\{p_j\}_{j=1}^\infty)}$ the $p$-shifted distribution on $S_n$ and $\text{proj}_n(\Pi)$ a $p$-shifted random permutation on $S_n$.
We note that in the case that $p_j=(1-q)q^{j-1}$, where $q\in(0,1)$, the measure $P^{(\{p_j\}_{j=1}^\infty)}_n$ is the Mallows distribution on $S_n$  with parameter $q$; see \cite{PT}, \cite{GO}.

\medskip

\noindent \bf Remark.\rm\ We assume in this paper that $p_j>0$, for all $j$. In fact, the $p$-shifted random permutation can
be constructed as long as $p_1>0$, with no positivity requirement on $p_j, j\ge2$. The positivity requirement
for all $j$ ensures that for all  $n$,  the support of the $p$-shifted measure $P_n$ is
all of $S_n$.
\medskip

It is known \cite{PT} that a random permutation under the $p$-shifted distribution  $P^{(\{p_j\}_{j=1}^\infty)}$  is \it  strictly regenerative\rm, where our
 definition of strictly regenerative is as follows.
For a permutation $\pi=\pi_{a+1}\pi_{a+2}\cdots\pi_{a+m}$, of $\{a+1,a+2,\cdots, a+m\}$, define
 $\text{red}(\pi)$, the reduced permutation of $\pi$, to be the permutation in $S_m$ given by
 $\text{red}(\pi)_i=\pi_{a+i}-m$.
 We will call a  random permutation $\Pi$ of $S_\infty$   \it strictly regenerative \rm if
 almost surely there exist
$0=T_0< T_1<T_2<\cdots$ such that $\Pi([T_j])=[T_j],\ j\ge1$,  and $\Pi([m])\neq[m]$ if $m\not\in\{T_1,T_2,\cdots\}$,  and such that the
random variables $\{T_k-T_{k-1}\}_{k=1}^\infty$ are IID and the  random permutations
$\{\text{red}(\Pi|_{[T_k]-[T_{k-1}]})\}_{k=1}^\infty$ are IID.
The numbers $\{T_n\}_{m=1}^\infty$ are called the renewal or regeneration numbers.
Our definition of strictly regenerative differs slightly from that in \cite{PT}.

Let $u_n$ denote the probability
that the $p$-shifted random permutation $\Pi$ has a renewal at the number $n$; that is,
 $u_n=P^{(\{p_j\}_{j=1}^\infty)}(\Pi([n])=[n])$.
 It follows easily from the construction of the random permutation that
\begin{equation}\label{un}
u_n=\prod_{j=1}^n(\sum_{i=1}^jp_i)=\prod_{j=1}^n(1-\sum_{i=j+1}^\infty p_i).
\end{equation}
See \cite{PT}.
Thus, $u_n>0$, for all $n$.
(Note that this
positivity, and the consequent aperiodicity of the renewal mechanism,
 does not require the positivity of all $p_j$, but only of $p_1$.)

The  strictly regenerative distribution $P^{(\{p_j\}_{j=1}^\infty)}$ is called \it positive recurrent\rm\ if $T_1$ has finite expectation:  $E^{(\{p_j\}_{j=1}^\infty)}\thinspace T_1<\infty$.
From standard renewal theory, it follows that
\begin{equation}\label{renewal}
\lim_{n\to\infty}u_n=\frac1{E^{(\{p_j\}_{j=1}^\infty)}\thinspace T_1}.
\end{equation}
Since $\sum_{j=1}^\infty\sum_{i=j+1}^\infty p_i=\sum_{j=1}^\infty jp_{j+1}$, it follows from \eqref{un} and \eqref{renewal} that
\begin{equation}\label{posreccond}
 P^{(\{p_j\}_{j=1}^\infty)} \ \text{is positive recurrent if and only if}\
\sum_{n=1}^\infty np_n<\infty.
\end{equation}

We  now state our theorem concerning super-clustering.

\begin{theorem}\label{supercluster}
Let $A^{(n)}_{l;k}\subset S_n$ denote the event that  the set of $l$ consecutive numbers $\{k,k+1,\cdots, k+l-1\}$
 appears in a set of $l$ consecutive positions. Let $\{p_n\}_{n=1}^\infty$ be a probability distribution on
 $\mathbb{N}$  with $p_j>0$, for all $j\in\mathbb{N}$.
 Also assume that the sequence $\{p_n\}_{n=1}^\infty$ is non-increasing.
Then for all $k\in\mathbb{N}$,
 \begin{equation}\label{fixedk}
\lim_{l\to\infty}\lim_{n\to\infty}P_n^{(\{p_j\}_{j=1}^\infty)}(A^{(n )}_{l,k})=\big(\prod_{j=1}^{k-1}\sum_{i=1}^jp_i\big)
 \big(\prod_{j=1}^\infty\sum_{i=1}^jp_i\big).
\end{equation}
Also, if $\lim_{n\to\infty}\min(k_n,n-k_n)=\infty$, then
  \begin{equation}\label{kn}
\lim_{l\to\infty}\liminf_{n\to\infty}P_n^{(\{p_j\}_{j=1}^\infty)}(A^{(n )}_{l,k_n})=
 \lim_{l\to\infty}\limsup_{n\to\infty}P_n^{(\{p_j\}_{j=1}^\infty)}(A^{(n )}_{l,k_n})=
 \big(\prod_{j=1}^\infty\sum_{i=1}^jp_i\big)^2.
\end{equation}
In particular, the limits in \eqref{fixedk} and \eqref{kn} are positive
if and only if $\sum_{n=1}^\infty np_n<\infty$, or equivalently, if and only if the $p$-shifted random permutation is positive recurrent.
\end{theorem}

\noindent \bf Remark.\rm\ If one removes the requirement that the sequence
$\{p_j\}_{j=1}^\infty$ be non-increasing, then it follows immediately from the proof of the theorem that
$$
\lim_{l\to\infty}\lim_{n\to\infty}P_n^{(\{p_j\}_{j=1}^\infty)}(A^{(n )}_{l,k})\ge\big(\prod_{j=1}^{k-1}\sum_{i=1}^jp_i\big)
 \big(\prod_{j=1}^\infty\sum_{i=1}^jp_i\big)
$$
and
$$
\lim_{l\to\infty}\liminf_{n\to\infty}P_n^{(\{p_j\}_{j=1}^\infty)}(A^{(n )}_{l,k_n})\ge
\big(\prod_{j=1}^\infty\sum_{i=1}^jp_i\big)^2.
$$
Thus, for this more general case, the finiteness of $\sum_{n=1}^\infty np_n$ is a sufficient condition for super-clustering.
\medskip

Consider Theorem \ref{supercluster} in the case of the Mallows
 distribution $P_n^q$ with parameter $q\in(0,1)$; that is, the case $p_j=(1-q)q^{j-1}$.
 From \eqref{duality},  \eqref{fixedk}
and \eqref{kn}, we have
\begin{equation}\label{3forMallows}
\begin{aligned}
&\lim_{l\to\infty}\lim_{n\to\infty}P_n^q(A^{(n )}_{l,k})=
\lim_{l\to\infty}\lim_{n\to\infty}P_n^q(A^{(n )}_{l,n+2-k-l})=\\
&\big(\prod_{j=1}^{k-1}(1-q^j)\big)\big(\prod_{j=1}^\infty(1-q^j)\big),\ \text{for all}\ k\in\mathbb{N};\\
&\lim_{l\to\infty}\lim_{n\to\infty}P_n^q(A^{(n )}_{l,k_n})=\big(\prod_{j=1}^\infty(1-q^j)\big)^2,
\text{if}\ \lim_{n\to\infty}\text{min}(k_n,n-k_n)=\infty.
\end{aligned}
\end{equation}
At the end of this section we prove the following easy asymptotic result.
\begin{proposition}\label{Cqasympt}
\begin{equation}\label{asymptCq}
\prod_{j=1}^\infty(1-q^j)\sim e^{-\frac{\pi^2}{6(1-q)}},\ \text{as}\ q\to1.
\end{equation}
\end{proposition}
 Proposition \ref{Cqasympt} and \eqref{3forMallows} yield the following corollary.
\begin{corollary}
If $\lim_{n\to\infty}\min(k_n,n-k_n)=\infty$, then
\begin{equation}\label{asymptCqprob}
\lim_{l\to\infty}\lim_{n\to\infty}P_n^q(A^{(n)}_{l;k_n})\sim e^{-\frac{\pi^2}{3(1-q)}}, \ \text{as}\ q\to1.
\end{equation}
\end{corollary}
\bf \noindent Remark.\rm\ In particular, if
$q_m=1-\frac c{\log m}$, then \eqref{asymptCqprob} gives
$$
\lim_{l\to\infty}\lim_{n\to\infty}P_n^{q_m}(A^{(n)}_{l;k_n})\sim m^{-\frac{\pi^2}{3c}},  \ \text{as}\ m\to\infty.
$$

\medskip

In \cite{P} we showed that under the $p$-shifted probability measure $P_n^{(\{p_j\}_{j=1}^\infty)}$,
the random variables $\{I_{<j}\}_{j=2}^\infty$ are independent,  and that  $1+I_{<k}$ is distributed as
 $\{p_j\}_{j=1}^\infty$, truncated at $k$:
 $$
 P_n^{(\{p_j\}_{j=1}^\infty)}(I_{<k}=i)=\frac{p_{i+1}}{\sum_{j=1}^kp_j},\
 \text{for}\ i=0,\cdots, k-1,\ \text{and}\   k=2,3,\cdots.
 $$
The statistics $\{I_{<k}\}_{k=2}^\infty$ are called the \it backward ranks\rm. As is well-known, a permutation is uniquely determined by its backward ranks.
This leads to an alternative way to construct a $p$-shifted random permutation in $S_n$ or in $S_\infty$.
Let $X$ be a random variable on $\mathbb{Z}^+$ whose distribution is characterized by $1+X$ having the distribution
$\{p_j\}_{j=1}^\infty$; that is,
\begin{equation}\label{X}
P(X=j)=p_{j+1},\ j=0,1,\cdots.
\end{equation}
Let $\{X_n\}_{n=2}^\infty$ be a sequence of independent random variables with the distribution of $X_n$ being the distribution of $X$ truncated at $n-1$:
\begin{equation}\label{Xn}
P(X_n=i)=\frac{p_{i+1}}{\sum_{j=1}^np_j},\ i=0,1,\cdots n-1.
\end{equation}
 To construct a $p$-shifted random permutation in $S_n$,  set
the number 1 down on a horizontal line. Now inductively, if the numbers  $\{1,\cdots, j-1\}$ have already been placed down on the line, where $2\le j\le n$,
then sample from $X_j$ independently of everything that has already occurred, and place
the number $j$ on the line in the position for which there are $X_j$ numbers to its right.
Thus, for example, to create a $p$-shifted random permutation in $S_4$,  if $X_2, X_3, X_4$ have been
sampled independently  as $X_2=1$, $X_3=2$ and $X_4=0$, then we obtain the permutation 3214.
To obtain a $p$-shifted random permutation in $S_\infty$, one just continues the above scenario indefinitely.
This alternative construction will be exploited for most of the proofs in this paper.

Since $EX=\sum_{j=1}^\infty jp_{j+1}$,
it follows from \eqref{posreccond} that
the $p$-shifted random permutation is positive recurrent  if and only if $EX<\infty$.
Note that for the random permutation  on $S_n$ or $S_\infty$ created in the previous paragraph, one has $X_j=I_{<j}$ for all appropriate $j$.
The total number of inversions in a permutation $\sigma\in S_n$ is given by $\mathcal{I}_n(\sigma):=\sum_{j=2}^nI_{<j}(\sigma)$. It follows from the construction in the above paragraph
that the inversion statistic $\mathcal{I}_n$ satisfies the following weak law of large numbers as $n\to\infty$:
\begin{equation}\label{wlln}
\frac{\mathcal{I}_n}n\ \text{under}\  P_n^{(\{p_j\}_{j=1}^\infty)}\ \text{converges in probability to}\ EX=\sum_{n=1}^\infty np_{n+1}\in(0,\infty].
\end{equation}

\medskip

\noindent \bf Remark 1.\rm\
In light of \eqref{wlln},
Theorem \ref{supercluster} shows that super-clustering occurs if and only if  the total inversion statistic
$\mathcal{I}_n$ has linear rather than super-linear growth.
\medskip

\bf\noindent Remark 2.\rm\  If $X^{(1)}$  and $X^{(2)}$ satisfy \eqref{X} with $\{p_j\}_{j=1}^\infty$ replaced respectively by distinct $\{p_j^{(1)}\}_{j=1}^\infty$ and $\{p_j^{(2)}\}_{j=1}^\infty$,
  and if $X^{(1)}$  stochastically
dominates $X^{(2)}$, that is, $\sum_{j=n}^\infty p_j^{(1)}\ge \sum_{j=n}^\infty p_j^{(2)}$, for all $n\in\mathbb{N}$, then it follows from
\eqref{fixedk} and \eqref{kn} that
the probability of super-clustering for the $p^{(2)}$-shifted random permutation is
 greater than for the  $p^{(1)}$-shifted random permutation.
This gives an explicit quantification of the inverse correlation between the tendency for inversion and the tendency for super-clustering.
\medskip

The considerations in this part of the paper lead naturally to the following characterization of the class of positive recurrent $p$-shifted distributions,
which might be of some independent interest.
\begin{proposition}\label{classif}
The class of $p$-shifted distributions, as $p$ runs over all probability distributions
$\{p_j\}_{j=1}^\infty$ whose  supports are all of $\mathbb{N}$, and that satisfy $\sum_{n=1}^\infty np_n<\infty$,
coincides with the class of probability distributions $P$ on $S_\infty$ that satisfy the following three conditions:

\noindent i. The backward ranks $\{I_{<j}\}_{j=2}^\infty$ are independent under $P$;

\noindent ii. A random permutation under $P$  is strictly regenerative with a positive recurrent
renewal mechanism, and the probability $u_1$ of renewal at the number 1 is positive;

\noindent iii. For all $n\in\mathbb{N}$, the support of $P_n(\cdot):=P\big(\text{\rm proj}_n^{-1}(\cdot)\big)$ is all of $S_n$.
\end{proposition}

\noindent\bf Remark.\rm\ The proof of the proposition also shows that if one removes the requirement that the support of the
distribution
$p$ is all of $\mathbb{N}$, and
only requires that $p_1>0$ (which in any case is necessary in order to implement  the $p$-shifted
construction), then the proposition holds with property (iii) deleted.

\medskip

We conclude this section with the proofs of Propositions \ref{dualprop} and \ref{Cqasympt}.

\noindent \it Proof of Proposition \ref{dualprop}.\rm\
We defined above the  reverse $\sigma^{\text{rev}}$
of a permutation $\sigma\in S_n$. The \it complement\rm\ of $\sigma$ is the permutation
$\sigma^{\text{com}}$ satisfying
$\sigma^{\text{com}}_i=n+1-\sigma_i,\ i=1,\cdots, n$.
Let $\sigma^{\text{rev-com}}$ denote the permutation obtained by applying reversal and then complementation to $\sigma$ (or equivalently, applying complentation and then reversal).
Since $\sigma^{\text{rev-com}}_i<\sigma^{\text{rev-com}}_j$ if and only $\sigma_{n+1-j}<\sigma_{n+1-i}$, it follows that $\sigma$ and $\sigma^{\text{rev-com}}$ have the same number of inversions, and thus, from the  definition of the Mallows distribution in \eqref{mallowsdef},
$P_n^q(\{\sigma\})=P_n^q(\{\sigma^{\text{rev-com}}\})$.
Using this along with the fact that $\sigma\in A^{(n)}_{l;k}$ if and only if $\sigma^{\text{rev-com}}\in A^{(n)}_{l;n+2-k-l}$
proves \eqref{duality}.\hfill $\square$
\medskip

\noindent \it Proof of Proposition \ref{Cqasympt}.\rm\
We have
\begin{equation}\label{Cq1}
\log\prod_{j=1}^\infty(1-q^j)=\sum_{j=1}^\infty\log(1-q^j),
\end{equation}
and
\begin{equation}\label{Cq2}
\int_1^\infty\log(1-q^x)dx\le\sum_{j=1}^\infty\log(1-q^j)\le\int_2^\infty\log(1-q^x)dx.
\end{equation}
Making the change of variables $y=q^x$ gives
\begin{equation}\label{Cq3}
\begin{aligned}
&\int_a^\infty\log(1-q^x)dx=-\frac1{\log q}\int_0^{q^a}\frac{\log(1-y)}ydy\sim\\
&\frac1{1-q}\big(\int_0^1\frac{\log(1-y)}ydy\big),\ \text{as}\ q\to1,\  \text{for}\ a>0.
\end{aligned}
\end{equation}
However, $\int_0^1\frac{\log(1-y)}ydy=-\int_0^1\big(\sum_{n=1}^\infty\frac{y^{n-1}}n\big)dy=-\sum_{n=1}^\infty\frac1{n^2}=
-\frac{\pi^2}6$.
Using this with \eqref{Cq1}-\eqref{Cq3}, we obtain \eqref{asymptCq}, proving the proposition. \hfill $\square$
\medskip

The alternative construction of $p$-shifted random permutations will be used for both the upper and lower bound calculations in the proof of Theorem \ref{supercluster}. The same type of upper bound calculations, specialized to the case of a Mallows distribution,   will also be used in the proofs of Theorems \ref{maintheorem} and  \ref{alpha1}. On the other hand, the original $p$-shifted construction, specialized to the case of a Mallows distribution, will be used for the lower bound calculations in the proof of Theorem \ref{maintheorem}.
In light of this, it will be convenient to begin with the proof of Theorem \ref{supercluster}, which is given in section \ref{super}. The proofs of Theorems \ref{maintheorem} and \ref{alpha1}  are given in sections
\ref{mainthmproof} and  \ref{alpha1proof}  respectively, and the proof of Proposition \ref{classif} is given
in section \ref{classification}.

\section{Proof of Theorem \ref{supercluster}}\label{super}
We note that the final statement of the theorem is almost immediate. Indeed,
$\sum_{i=1}^jp_i=1-\sum_{i=j+1}^\infty p_i$ and $\sum_{j=1}^\infty\big(\sum_{i=j+1}^\infty p_i\big)=\sum_{j=1}^\infty jp_{j+1}$.

We now turn to the proofs of \eqref{fixedk} and \eqref{kn}.
We use the alternative method for constructing the $p$-shifted random permutation, as described after
\eqref{X}.
Thus, we consider a sequence of independent random variables $\{X_n\}_{n=2}^\infty$, with
$X_n$ distributed as in \eqref{Xn}.
For the proof, we will use the notation
\begin{equation}\label{Nn}
N_n=\sum_{i=1}^np_i=P(X\le n-1),\  n\in\mathbb{N},\ \text{and}\ \ N_0=0,
\end{equation}
where $X$ is as in \eqref{X}.
Note that $N_n$ is the normalization constant on the right hand side of \eqref{Xn}.
Although $P_n^{(\{p_j\}_{j=1}^\infty)}$ denotes the $p$-shifted probability measure on $S_n$, we will also use this notation for  probabilities of events related to the random variables
$\{X_j\}_{j=2}^n$. However, probabilities of events related to $X$ will still be denoted by $P$.

We begin with the proof of \eqref{fixedk}. Fix $k\in\mathbb{N}$.
Consider the event, which we denote by $B_{l;k}$,  that after the
first $k+l-1$ positive integers have been placed down on the horizontal line,
the set  of $l$ numbers $\{k,k+1,\cdots, k+l-1\}$ appear in a set of $l$ consecutive positions.
Then $B_{l;k}=\cup_{a=0}^{k-1}B_{l;k;a}$, where the events $\{B_{l;k;a}\}_{a=0}^{k-1}$
are disjoint, with $B_{l;k;a}$ being the event that the set  of $l$ numbers $\{k,k+1,\cdots, k+l-1\}$
appear in a set of $l$ consecutive positions and also that exactly $a$ of the numbers in $[k-1]$ are to the right of this set.
We calculate $P_n^{(\{p_j\}_{j=1}^\infty)}(B_{l;k;a})$.

Suppose that we have already placed  down on the horizontal line the numbers in $[k-1]$.
Their relative positions are irrelevant for our considerations.
Now we use $X_k$ to insert  on the line the number $k$. Suppose that $X_k=a$, $a\in\{0,\cdots, k-1\}$.
Then the number $k$ is inserted  on the line in the position for which  $a$ of the numbers in  $[k-1]$  are to its right. Now in order for $k+1$ to be placed in a position adjacent to $k$, we need
$X_{k+1}\in\{a,a+1\}$. (If $X_{k+1}=a$, then $k+1$ will appear directly to the right of $k$, while if $X_{k+1}=a+1$, then $k+1$ will appear directly to the left of $k$.)
If this occurs, then $\{k,k+1\}$ are adjacent, and $a$ of the  numbers in $[k-1]$ are to the right of $\{k,k+1\}$.
Continuing in this vein, for $i\in\{1,\cdots, l-2\}$,
given that the numbers $\{k,\cdots, k+i\}$  are adjacent to one another, and $a$ of the numbers in $[k-1]$ appear
to the right of $\{k,\cdots, k+i\}$, then in order for
$k+i+1$ to be placed
so that $\{k,\cdots, k+i+1\}$ are all adjacent to one another (with $a$ of the numbers in $[k-1]$ appearing to the right of these numbers), we need $X_{k+i+1}\in\{a,\cdots, a+i+1\}$.
We conclude then that
$P_n^{(\{p_j\}_{j=1}^\infty)}(B_{l;k;a})=\prod_{j=0}^{l-1}P_n^{(\{p_j\}_{j=1}^\infty)}(X_{k+j}\in\{a,\cdots, a+j\})$.
Using \eqref{Xn},
we have
\begin{equation}\label{Ba}
P_n^{(\{p_j\}_{j=1}^\infty)}(B_{l;k;a})=\prod_{j=0}^{l-1}P_n^{(\{p_j\}_{j=1}^\infty)}(X_{k+j}\in\{a,\cdots, a+j\})=
\prod_{j=0}^{l-1}\frac{N_{a+j+1}-N_a}{N_{k+j}}.
\end{equation}

We now consider the conditional probability,
$P_n^{(\{p_j\}_{j=1}^\infty)}(A^{(n)}_{l;k}|B_{l;k;a})$, that is, the probability,
 given that
 $B_{l;k;a}$ has occurred,  that the numbers $k+l,\cdots, n$ are inserted in such a way so as to preserve the mutual adjacency of the numbers in the set
$\{k,\cdots,k+l-1\}$.
We will obtain lower and upper bounds on this conditional probability. However, first we note that it is clear from the  construction that
$P_n^{(\{p_j\}_{j=1}^\infty)}(A^{(n)}_{l;k}|B_{l;k;a})$ is decreasing in $n$. Thus, since
$P_n^{(\{p_j\}_{j=1}^\infty)}(B_{l;k;a})$ is independent of $n$,
it follows that
$P_n^{(\{p_j\}_{j=1}^\infty)}(A^{(n)}_{l;k})$ is decreasing in $n$. Consequently $\lim_{n\to\infty}P_n^{(\{p_j\}_{j=1}^\infty)}(A^{(n)}_{l;k})$ exists.

We now turn to a lower bound on $P_n^{(\{p_j\}_{j=1}^\infty)}(A^{(n)}_{l;k}|B_{l;k;a})$.
Our lower bound will be the probability of the event  that all of the remaining numbers are inserted to the right of the set $\{k,\cdots,k+l-1\}$.
This event is given by $\cap_{j=0}^{n-k-l}\{X_{k+l+j}\le a+j\}$.
Thus, we have
\begin{equation}\label{firstlower}
\begin{aligned}
P_n^{(\{p_j\}_{j=1}^\infty)}(A^{(n)}_{l;k}|B_{l;k;a})\ge P_n^{(\{p_j\}_{j=1}^\infty)}(\cap_{j=0}^{n-k-l}\{X_{k+l+j}\le a+j\})=
\prod_{j=0}^{n-k-l}\frac{N_{a+j+1}}{N_{k+l+j}}.
\end{aligned}
\end{equation}

Writing $P_n^{(\{p_j\}_{j=1}^\infty)}(A^{(n)}_{l;k})=\sum_{a=0}^{k-1}P_n^{(\{p_j\}_{j=1}^\infty)}(B_{l;k;a})
P_n^{(\{p_j\}_{j=1}^\infty)}(A^{(n)}_{l;k}|B_{l;k;a})$,
  \eqref{Ba} and \eqref{firstlower} yield
\begin{equation}\label{firstlowerbound}
P_n^{(\{p_j\}_{j=1}^\infty)}(A^{(n)}_{l;k})\ge\sum_{a=0}^{k-1}\big(\prod_{j=0}^{l-1}\frac{N_{a+j+1}-N_a}{N_{k+j}}\big)
\big(\prod_{j=0}^{n-k-l}\frac{N_{a+j+1}}{N_{k+l+j}}\big).
\end{equation}
We have $\prod_{j=0}^{n-k-l}\frac{N_{a+j+1}}{N_{k+l+j}}=\frac{N_{a+1}\cdots n_{k+l-1}}{N_{n-k-l+a+2}\cdots N_n}$.
Using this along with the fact that $\lim_{n\to\infty}N_n=1$ and
the fact that the limit on the left hand side of
\eqref{firstlowerbound} exists, we have
\begin{equation}\label{firstlowerboundinfinity}
\begin{aligned}
&\lim_{n\to\infty}P_n^{(\{p_j\}_{j=1}^\infty)}(A^{(n)}_{l;k})\ge
\sum_{a=0}^{k-1}\big(\prod_{j=0}^{l-1}\frac{N_{a+j+1}-N_a}{N_{k+j}}\big)\big(\prod_{i=a+1}^{k+l-1}N_i\big)=\\
&\sum_{a=0}^{k-1}\big(\prod_{j=0}^{l-1}(N_{a+j+1}-N_a)\big)\big(\prod_{i=a+1}^{k-1}N_i\big).
\end{aligned}
\end{equation}
We now let $l\to\infty$ in \eqref{firstlowerboundinfinity}.
We only consider the term in the summation with $a=0$, because
it turns out that the terms  with $a\ge1$
converge to 0 as $l\to\infty$.
We obtain
\begin{equation}\label{lowerboundthm1}
\lim_{l\to\infty}\lim_{n\to\infty}P_n^{(\{p_j\}_{j=1}^\infty)}(A^{(n)}_{l;k})\ge \big(\prod_{j=1}^{k-1}N_j\big)\big(\prod_{j=1}^\infty N_j\big)=
\big(\prod_{j=1}^{k-1}\sum_{i=1}^j p_i\big)\big(\prod_{j=1}^\infty\sum_{i=1}^jp_i\big).
\end{equation}

Note that by the assumption that $\{p_n\}_{n=1}^\infty$ is non-increasing,
it follows that $P(X\not\in\{j+1,\cdots,j+l-1\})$ is increasing in $j$.
 Also, note that $P(X\not\in\{j+1,\cdots,j+l-1\})>P_n^{(\{p_j\}_{j=1}^\infty)}(X_m\not\in\{j+1,\cdots,j+l-1\})$, for $j+l\le m$.
These facts will be used as we turn now to
 an upper bound on $P_n^{(\{p_j\}_{j=1}^\infty)}(A^{(n)}_{l;k}|B_{l;k;a})$,
the conditional probability given $B_{l;k;a}$ that the numbers $k+l,\cdots, n$ are inserted in such a way so as to
preserve the mutual adjacency of the set
$\{k,\cdots,k+l-1\}$.
 First the number $k+l$ is inserted. The probability that its insertion preserves the mutual adjacency property of the set
$\{k,\cdots,k+l-1\}$
is $P_n^{(\{p_j\}_{j=1}^\infty)}(X_{k+l}\not\in\{a+1,\cdots, a+l-1\})$, which is less than
$P(X\not\in\{a+1,\cdots, a+l-1\})$. If the insertion of $k+l$ preserves the mutual adjacency, then either
$X_{k+l}\in\{0,\cdots, a\}$ or $X_{k+l}\in\{a+l,\cdots, k+l-1\}$.
If $X_{k+l}\in\{0,\cdots, a\}$, then in order for the mutually adjacency to be preserved when the number
$k+l+1$ is inserted, one needs $\big\{X_{k+1+1}\not\in\{a+2,\cdots,a+l\}\big\}$, while if
$X_{k+l}\in\{a+l,\cdots, k+l-1\}$, then one needs
$\big\{X_{k+1+1}\not\in\{a+1,\cdots,a+l-1\}\big\}$.
The probability of either of these events is less than
$P(X\not\in\{a+2,\cdots, a+l\})$. Thus, an upper bound for the conditional probability given $B_{l;k;a}$ that
the insertion of $k+l$ and $k+l+1$ preserves the mutual adjacency is
$P(X\not\in\{a+1,\cdots, a+l-1\})P(X\not\in\{a+2,\cdots, a+l\})$.
Continuing in this vein, we conclude that
\begin{equation}\label{firstupper}
\begin{aligned}
&P_n^{(\{p_j\}_{j=1}^\infty)}(A^{(n)}_{l;k}|B_{l;k;a})\le \prod_{j=1}^{n-k-l+1} P(X\not\in\{a+j,\cdots, a+j+l-2\})=\\
&\prod_{j=1}^{n-k-l+1}(1-N_{a+j+l-1}+N_{a+j}).
\end{aligned}
\end{equation}

From \eqref{Ba}  and \eqref{firstupper},  we have
\begin{equation}\label{firstupperbound}
P_n^{(\{p_j\}_{j=1}^\infty)}(A^{(n)}_{l;k})\le\sum_{a=0}^{k-1}\big(\prod_{j=0}^{l-1}\frac{N_{a+j+1}-N_a}{N_{k+j}}\big)
\big(\prod_{j=1}^{n-k-l+1}(1-N_{a+j+l-1}+N_{a+j})\big).
\end{equation}
Letting $n\to\infty$ and using the fact that the limit on the left hand side exists, we have
 \begin{equation}\label{another}
\lim_{n\to\infty}P_n^{(\{p_j\}_{j=1}^\infty)}(A^{(n)}_{l;k})\le
\sum_{a=0}^{k-1}\big(\prod_{j=0}^{l-1}\frac{N_{a+j+1}-N_a}{N_{k+j}}\big)
\big(\prod_{j=1}^\infty(1-N_{a+j+l-1}+N_{a+j})\big).
\end{equation}
For $a\in\{1,\cdots, k-1\}$, we have  $\frac{N_{a+j+1}-N_a}{N_{k+j}}<1-N_a\in(0,1)$, for all $j\ge0$.
Therefore, when letting $l\to\infty$ in \eqref{another}, a contribution will come from the right hand side only when $a=0$. We obtain
\begin{equation}\label{upperboundthm1}
\begin{aligned}
&\lim_{l\to\infty}\lim_{n\to\infty}P_n^{(\{p_j\}_{j=1}^\infty)}(A^{(n)}_{l;k})\le\lim_{l\to\infty}
\big(\prod_{j=0}^{l-1}\frac{N_{j+1}}{N_{k+j}}\big)
\big(\prod_{j=1}^\infty(1-N_{j+l-1}+N_j)\big)=\\
& \big(\prod_{j=1}^{k-1}N_j\big)\big(\prod_{j=1}^\infty N_j\big)=
\big(\prod_{j=1}^{k-1}\sum_{i=1}^j p_i\big)\big(\prod_{j=1}^\infty\sum_{i=1}^jp_i\big).
\end{aligned}
\end{equation}
Now \eqref{fixedk} follows from \eqref{lowerboundthm1} and \eqref{upperboundthm1}.

We now turn to the proof of \eqref{kn}.  As with the proof of \eqref{fixedk}, the term with $a=0$ will dominate.
Thus, for the lower bound, using \eqref{Ba} and \eqref{firstlower} with $k=k_n$ and ignoring the terms with $a\ge1$, we have
\begin{equation}\label{yetanother}
\begin{aligned}
&P_n^{(\{p_j\}_{j=1}^\infty)}(A^{(n)}_{l;k_n})\ge\big(\prod_{j=0}^{l-1}\frac{N_{j+1}}{N_{k_n+j}}\big)
\big(\prod_{j=0}^{n-k_n-l}\frac{N_{j+1}}{N_{k_n+l+j}}\big).
\end{aligned}
\end{equation}
Letting $n\to\infty$ in \eqref{yetanother} and
  using the assumption that $\lim_{n\to\infty}\min(k_n,n-k_n)=\infty$, it follows that
$$
\liminf_{n\to\infty}P_n^{(\{p_j\}_{j=1}^\infty)}(A^{(n)}_{l;k_n})\ge\big(\prod_{j=1}^lN_j\big)\big(\prod_{j=1}^\infty N_j\big).
$$
Now letting $l\to\infty$ gives
\begin{equation}\label{lowerboundthm1infty}
\lim_{l\to\infty}\liminf_{n\to\infty}P_n^{(\{p_j\}_{j=1}^\infty)}(A^{(n)}_{l;k_n})\ge \big(\prod_{j=1}^\infty N_j\big)^2=
\big(\prod_{j=1}^\infty\sum_{i=1}^jp_i\big)^2.
\end{equation}

For the upper bound, let $k=k_n$ in
\eqref{firstupperbound}.
The second factor in the summand
$\big(\prod_{j=0}^{l-1}\frac{N_{a+j+1}-N_a}{N_{k_n+j}}\big)
\big(\prod_{j=1}^{n-k_n-l+1}(1-N_{a+j+l-1}+N_{a+j})\big)$
is less than 1, while the first factor  in the summand satisfies
$$
\prod_{j=0}^{l-1}\frac{N_{a+j+1}-N_a}{N_{k_n+j}}\le\frac{N_{a+1}-N_a}{N_{k_n}}=\frac{p_{a+1}}{N_{k_n}}\le\frac{p_{a+1}}{p_1},
$$
 for $a\in\{0,\cdots, k_n-1\}$ and $n\ge1$.
Since $\sum_{a=0}^\infty\frac{p_{a+1}}{p_1}<\infty$,  the dominated convergence theorem and the assumption that
 $\lim_{n\to\infty}\min(k_n,n-k_n)=\infty$
 allow us to conclude upon letting $n\to\infty$ in \eqref{firstupperbound} with $k=k_n$ that
 \begin{equation}\label{almostfinaltime}
\limsup_{n\to\infty}P_n^{(\{p_j\}_{j=1}^\infty)}(A^{(n)}_{l;k_n})\le\sum_{a=0}^\infty\big(\prod_{j=0}^{l-1}(N_{a+j+1}-N_a)\big)\big(\prod_{j=1}^\infty(1-N_{a+j+l-1}+N_{a+j})\big).
\end{equation}
For $a\ge1$, we have $N_{a+j+1}-N_a\in(0,1-p_1)$. Consequently, when letting $l\to\infty$ in \eqref{almostfinaltime},
 a contribution will come from the right hand side only when $a=0$. We obtain
\begin{equation}\label{upperboundthm1infty}
\lim_{l\to\infty}\limsup_{n\to\infty}P_n^{(\{p_j\}_{j=1}^\infty)}(A^{(n)}_{l;k_n})\le\big(\prod_{j=1}^\infty N_j\big)^2=
\big(\prod_{j=1}^\infty\sum_{i=1}^jp_i\big)^2.
\end{equation}
Now \eqref{kn}  follows from \eqref{lowerboundthm1infty} and  \eqref{upperboundthm1infty}.
\hfill $\square$

\section{Proof  of Theorem \ref{maintheorem}}\label{mainthmproof}
We will prove \eqref{mainresultgen} and \eqref{mainresult} in tandem. Note that the lower bounds in \eqref{mainresultgen} and \eqref{mainresult} are the same; only the upper bounds differ.
Recall  that \eqref{mainresult} is stated to hold under the assumption
$\frac{\min(k_n,n-k_n)}{n^\alpha}=\infty$, while  \eqref{mainresultgen} is stated to hold with no assumption on $\{k_n\}_{n=1}^\infty$.
Thus, we need to prove the common lower bound in
\eqref{mainresultgen} and \eqref{mainresult},
 as well as the upper bound in \eqref{mainresultgen},   with no assumption on $\{k_n\}_{n=1}^\infty$, while we need to prove the upper bound in
 \eqref{mainresult} under the above noted assumption
on $\{k_n\}_{n=1}^\infty$.
In fact, for our proofs, we will always need   to assume that
\begin{equation}\label{knassump}
\lim_{n\to\infty}\frac{k_n}{n^\alpha}=\infty.
\end{equation}
What allows us to make this assumption is Proposition \ref{dualprop}.
Thus,  in the sequel we will always assume that \eqref{knassump} holds.

For the upper bound, we follow the same construction used in the upper bound in Theorem \ref{supercluster}.
We start from \eqref{firstupperbound} with $k$ and $q$ replaced by $k_n$ and $q_n$. Since the Mallows distribution with parameter $q_n$ is the $p$-shifted distribution with
$p_j=(1-q_n)q_n^{j-1}$, it follows from \eqref{Nn} that for the case at hand,
\begin{equation}\label{Nb}
N_b=\sum_{i=1}^b(1-q_n)q_n^{i-1}=1-q_n^b.
\end{equation}
Substituting \eqref{Nb} in \eqref{firstupperbound}, we obtain
\begin{equation}\label{firstupperboundagain}
P_n^{q_n}(A^{(n)}_{l;k_n})\le\prod_{j=0}^{l-1}\frac{1-q_n^{j+1}}{1-q_n^{k_n+j}}\sum_{a=0}^{k_n-1}q_n^{al}
\prod_{j=1}^{n-k_n-l+1}\big(1-q_n^{a+j}+q_n^{a+j+l-1}\big).
\end{equation}
We have
\begin{equation}\label{qntob}
1-q_n^b=1-(1-\frac c{n^\alpha})^b\sim \frac{bc}{n^\alpha}, \text{for}\ b\in \mathbb{N},
\end{equation}
and
\begin{equation}\label{qnkn}
1-q_n^{k_n+j}=1-(1-\frac c{n^\alpha})^{k_n+j}\ge1-e^{-\frac{c(k_n+j)}{n^\alpha}}.
\end{equation}
From \eqref{qntob} and \eqref{qnkn} along with the assumption on  $q_n$ and the assumption \eqref{knassump} on $k_n$, the term multiplying the summation in \eqref{firstupperboundagain} satisfies
\begin{equation}\label{multsum}
\prod_{j=0}^{l-1}\frac{1-q_n^{j+1}}{1-q_n^{k_n+j}}\sim \frac{l!c^l}{n^{\alpha l}}.
\end{equation}
Using \eqref{qntob}, the summation in \eqref{firstupperboundagain} satisfies
\begin{equation}\label{summationterm}
\sum_{a=0}^{k_n-1}
q_n^{al}\prod_{j=1}^{n-k_n-l+1}\big(1-q_n^{a+j}+q_n^{a+j+l-1}\big)\sim
\sum_{a=0}^{k_n-1}
q_n^{al}\prod_{j=1}^{n-k_n-l+1}\Big(1-\frac{q_n^{a+j}(l-1)c}{n^\alpha}\Big).
\end{equation}

We split up the continuation of the proof of the upper bound between the case that no assumption is made on $k_n$
(accept for \eqref{knassump}, as explained above), in which case  we need to prove the upper bound in \eqref{mainresultgen},
and the case that  $k_n$ is assumed to satisfy $\frac{\min(k_n,n-k_n)}{n^\alpha}=\infty$, in which case we need to prove the upper bound in \eqref{mainresult}.
We begin with the former case.
In this case, from \eqref{firstupperboundagain}  along with \eqref{qntob}, \eqref{multsum} and \eqref{summationterm}, we have
\begin{equation*}
P_n^q(A^{(n)}_{l;k_n})\lesssim\frac{l!c^l}{n^{\alpha l}}\sum_{a=0}^{k_n-1}q_n^{al}\le
\frac{l!c^l}{n^{\alpha l}}\frac1{1-q_n^{l}}\sim\frac1l\thinspace\frac{l!c^{l-1}}{n^{\alpha(l-1)}},
\end{equation*}
which is the upper bound in \eqref{mainresultgen}.

Now consider the case  that $k_n$ is assumed to satisfy $\frac{\min(k_n,n-k_n)}{n^\alpha}=\infty$, in which case we need to prove the upper bound in \eqref{mainresult}. In the previous case, we simply replaced the product on the right hand side of
\eqref{summationterm} by one.
For the current case, we analyze this product.
We write
\begin{equation}\label{loganalysis}
\log\prod_{j=1}^{n-k_n-l+1}\Big(1-\frac{q_n^{a+j}(l-1)c}{n^\alpha}\Big)=\sum_{j=1}^{n-k_n-l+1}
\log\Big(1-\frac{q_n^{a+j}(l-1)c}{n^\alpha}\Big).
\end{equation}
We have
\begin{equation}\label{logint}
\begin{aligned}
&\int_0^{n-k_n-l+1} \log\Big(1-\frac{q_n^{a+x}(l-1)c}{n^\alpha}\Big)dx
\le \sum_{i=1}^{n-k_n-l+1}\log\Big(1-\frac{q_n^{a+i}(l-1)c}{n^\alpha}\Big)\le\\
&\int_1^{n-k_n-l+2}\log\Big(1-\frac{q_n^{a+x}(l-1)c}{n^\alpha}\Big)dx.
\end{aligned}
\end{equation}
Making the change of variables, $y=q_n^x$, we have
\begin{equation}\label{intAB}
\int_A^B\log\Big(1-\frac{q_n^{a+x}(l-1)c}{n^\alpha}\Big)dx=
-\frac1{\log q_n}\int_{q_n^B}^{q_n^A}\frac{\log\big(1-\frac{q_n^a(l-1)c}{n^\alpha}y\big)}ydy.
\end{equation}
From \eqref{intAB} and the assumptions on $q_n$ and  $k_n$, both the left and the right hand sides of \eqref{logint} are asymptotic
to
$\frac{n^\alpha}c\int_0^1\frac{\log\big(1-\frac{q_n^a(l-1)c}{n^\alpha}y\big)}ydy$, which in turn is asymptotic to
$-(l-1)q_n^a$, uniformly over $a\in\{0,\cdots, k_n-1\}$.
Using this with  \eqref{loganalysis} and \eqref{logint} gives
\begin{equation}\label{productasymp}
\prod_{j=1}^{n-k_n-l+1}\Big(1-\frac{q_n^{a+j}(l-1)c}{n^\alpha}\Big)\sim e^{-(l-1)q_n^a},\ \text{uniformly over}\  a\in\{0,\cdots, k_n-1\}.
\end{equation}

From \eqref{firstupperboundagain}  along with \eqref{multsum}, \eqref{summationterm} and \eqref{productasymp}, we obtain
\begin{equation}\label{innermainupperest}
P_n^q(A^{(n)}_{l;k_n})\lesssim\frac{l!c^l}{n^{\alpha l}}\sum_{a=0}^{k_n-1}q_n^{al}e^{-(l-1)q_n^a}.
\end{equation}
By the assumptions on $k_n$ and $q_n$, $\sum_{a=0}^{k_n-1}q_n^{al}e^{-(l-1)q_n^a}$ is asymptotic to
\newline
$\int_0^{k_n}q_n^{xl}e^{-(l-1)q_n^x}dx$. Making the change of variables $y=q_n^x$,
this integral is equal to $-\frac1{\log q_n}\int_{q_n^{k_n}}^1y^{l-1}e^{-(l-1)y}dy$, which in turn is asymptotic to
\newline $\frac{n^\alpha}c\int_0^1y^{l-1}e^{-(l-1)y}dy$. Thus,
\begin{equation}\label{sumasymptotics}
\sum_{a=0}^{k_n-1}q_n^{al}e^{-(l-1)q_n^a}\sim \frac{n^\alpha}c\int_0^1y^{l-1}e^{-(l-1)y}dy.
\end{equation}
From \eqref{innermainupperest} and \eqref{sumasymptotics}, we conclude that
\begin{equation*}
P_n^q(A^{(n)}_{l;k_n})\lesssim \big(\int_0^1y^{l-1}e^{-(l-1)y}dy\big)\thinspace\frac{c^{l-1}l!}{n^{\alpha(l-1)}},
\end{equation*}
which is the upper bound in \eqref{mainresult}.

We now turn to the lower bound.
Our only assumption on $k_n$ is \eqref{knassump}.
The method used
in the proof of  Theorem \ref{supercluster} and in the proof of the upper bound here,
via the alternative method for constructing a $p$-shifted random permutation,
is not precise enough to be of use in the proof of
the lower bound here. For this, we utilize the original construction for $p$-shifted random permutations on $S_n$, specializing to the   Mallows distribution with parameter $q_n$, for which
$p_j=(1-q_n)q_n^{j-1}$. We   use the notation  $P_n^{q_n}$ not only for the Mallows distribution itself, but also for probabilities of events associated with
 the  construction.
With regard to this construction, for $j\in\{0,\cdots, k_n-1\}$,  let $C_{j;k_n,l}$  denote the event that exactly $j$ numbers from the set $\{1,\cdots, k_n-1\}$ appear in the permutation before any number from the set
$\{k_n,\cdots, k_n+l-1\}$ appears.
We  calculate $P_n^{q_n}(C_{j;k_n,l})$ explicitly.
For $a,b\in\mathbb{N}$, let $r_{a,b}$
denote the probability that in the  construction,
the first number that appears from the set $\{1,\cdots, a+b\}$ comes from the set $\{1,\cdots, a\}$.
Then
\begin{equation}\label{rab}
r_{a,b}=\frac{\sum_{j=1}^a(1-q_n)q_n^{j-1}}{\sum_{j=1}^{a+b}(1-q_n)q_n^{j-1}}=\frac{1-q_n^a}{1-q_n^{a+b}}.
\end{equation}
For convenience, define $r_{0,b}=0$. Then from the construction, it follows that
\begin{equation}\label{Cj}
P_n^{q_n}(C_{j;k_n,l})=\big(\prod_{i=1}^jr_{k_n-i,l}\big)(1-r_{k_n-j-1,l}),\ j=0,\cdots, k_n-1.
\end{equation}
From \eqref{rab} and \eqref{Cj}, we have
\begin{equation}\label{Cjexplicit}
\begin{aligned}
&P_n^{q_n}(C_{j;k_n,l})=\big(\prod_{i=1}^j\frac{1-q_n^{k_n-i}}{1-q_n^{k_n-i+l}}\big)\frac{q_n^{k_n-j-1}-q_n^{k_n-j-1+l}}{1-q_n^{k_n-j-1+l}}=\\
&\frac{(1-q_n^l)q_n^{k_n-1-j}}{1-q_n^{k_n-1-j+l}}\thinspace \frac{\prod_{b=k_n-j}^{\min(k_n-j+l-1,k_n-1)}(1-q_n^b)}{\prod_{b=\max(k_n-j+l,k_n)}^{k_n+l-1}(1-q_n^b)}=\\
&\begin{cases}\frac{(1-q_n^l)q_n^{k_n-1-j}}{1-q_n^{k_n-1-j+l}}\thinspace \frac{\prod_{b=k_n-j}^{k_n-1}(1-q_n^b)}{\prod_{b=k_n-j+l}^{k_n+l-1}(1-q_n^b)}, \ j\le l-1;\\
\frac{(1-q_n^l)q_n^{k_n-1-j}}{1-q_n^{k_n-1-j+l}}\thinspace \frac{\prod_{b=k_n-j}^{k_n-j+l-1}(1-q_n^b)}{\prod_{b=k_n}^{k_n+l-1}(1-q_n^b)}, \ j\ge l,\end{cases}
\ j=0,\cdots, k_n-1.
\end{aligned}
\end{equation}

In order for the event $A_{l;k_n}^{(n)}$  to occur, the $l$ numbers $\{k_n,\cdots, k_n+l-1\}$
 must appear consecutively (in arbitrary order) in the construction. Thus,
given  the event $C_{j;k_n,l}$,  in order for the event $A_{l;k_n}^{(n)}$
to occur, all of the other $l-1$ numbers in $\{k_n,\cdots, k_n+l-1\}$ must occur immediately after the appearance of the first number from this set.
Given $C_{j;k_n,l}$, after the appearance of the first number from $\{k_n,\cdots, k_n+l-1\}$, there are still $k_n-1-j$ numbers from $\{1,\cdots, k_n-1\}$ that have not yet appeared, as well as
a certain amount of numbers from $\{k_n+l,\cdots, n\}$.
Thus, a lower bound on $P_n^{q_n}(A_{l;k_n}^{(n)}|C_{j;k_n,l})$ is obtained by assuming that
none of the numbers from $\{k_n+l,\cdots, n\}$ have  yet appeared.
(Here it is appropriate to note that if we calculate an upper bound by assuming that all of the numbers from
 $\{k_n+l,\cdots, n\}$ have  already appeared,  then the upper bound we arrive at for $P_n^{q_n}(A_{l;k_n}^{(n)})$ is not as good as
the upper bound in \eqref{mainresult}.)

In order to calculate explicitly this lower bound,
for $a,b,c\in\mathbb{N}$, let $r_{a,b,c}$ denote the probability that
the first number that appears from the set $\{1,\cdots, a+b+c\}$ comes from the set
$\{1,\cdots, a\}\cup\{a+b+1,\cdots, a+b+c\}$. Then
$$
r_{a,b,c}=\frac{\sum_{j=1}^a(1-q_n)q_n^{j-1}+\sum_{j=a+b+1}^{a+b+c}(1-q_n)q_n^{j-1}}{\sum_{j=1}^{a+b+c}(1-q_n)q_n^{j-1}}=\frac{1-q_n^a+q_n^{a+b}-q_n^{a+b+c}}{1-q_n^{a+b+c}}.
$$
From the construction, the lower bound on $P_n^{q_n}(A_{l;k_n}^{(n)}|C_{j;k_n,l})$,  obtained by assuming that
none of the numbers from $\{k_n+l,\cdots, n\}$ have  yet appeared, is given by
\begin{equation}\label{innerlowerboundmain}
\begin{aligned}
&P_n^{q_n}(A_{l;k_n}^{(n)}|C_{j;k_n,l})\ge\prod_{i=1}^{l-1}(1-r_{k_n-1-j,i,n-k_n-l+1})=\prod_{i=1}^{l-1}\frac{q_n^{k_n-1-j}-q_n^{k_n-1-j+i}}{1-q_n^{n-l-j+i}}=\\
&q_n^{(l-1)(k_n-1-j)}\frac{\prod_{b=1}^{l-1}(1-q_n^b)}{\prod_{b=n-l-j+1}^{n-j-1}(1-q_n^b)}.
\end{aligned}
\end{equation}

From \eqref{Cjexplicit} and \eqref{innerlowerboundmain},  we have
\begin{equation}\label{totallowerbound}
\begin{aligned}
&P_n^{q_n}(A_{l;k_n}^{(n)})=\sum_{j=0}^{k_n-1}P_n^{q_n}(C_{j;k_n,l})P_n^{q_n}(A_{l;k_n}^{(n)}|C_{j;k_n,l})\ge\\
&\sum_{j=l}^{k_n-1} \frac{(1-q_n^l)q_n^{k_n-1-j}}{1-q_n^{k_n-1-j+l}}\thinspace \frac{\prod_{b=k_n-j}^{k_n-j+l-1}(1-q_n^b)}{\prod_{b=k_n}^{k_n+l-1}(1-q_n^b)}
q_n^{(l-1)(k_n-1-j)}\frac{\prod_{b=1}^{l-1}(1-q_n^b)}{\prod_{b=n-l-j+1}^{n-j-1}(1-q_n^b)}.
\end{aligned}
\end{equation}
By the assumption on  $q_n$, the right hand side of \eqref{totallowerbound} satisfies
\begin{equation}\label{asymprhs}
\begin{aligned}
&\sum_{j=l}^{k_n-1} \frac{(1-q_n^l)q_n^{k_n-1-j}}{1-q_n^{k_n-1-j+l}}\thinspace \frac{\prod_{b=k_n-j}^{k_n-j+l-1}(1-q_n^b)}{\prod_{b=k_n}^{k_n+l-1}(1-q_n^b)}
q_n^{(l-1)(k_n-1-j)}\frac{\prod_{b=1}^{l-1}(1-q_n^b)}{\prod_{b=n-l-j+1}^{n-j-1}(1-q_n^b)}\gtrsim\\
&\prod_{b=1}^l(1-q_n^b)\sum_{j=l}^{k_n-1}q_n^{l(k_n-1-j)}\prod_{b=k_n-j}^{k_n-j+l-2}(1-q_n^b)\sim\frac{l!c^l}{n^{l\alpha}}\sum_{j=l}^{k_n-1}q_n^{l(k_n-1-j)}\big(1-q_n^{k_n-j}\big)^{l-1}.
\end{aligned}
\end{equation}
And
\begin{equation}
\sum_{j=l}^{k_n-1}q_n^{l(k_n-1-j)}\big(1-q_n^{k_n-j}\big)^{l-1}\sim\int_0^{k_n-1-l} q_n^{xl}(1-q_n^x)^{l-1}dx.
\end{equation}
Making the change of variables $y=q_n^x$, and using the assumption on $q_n$ and the assumption on $k_n$ in
\eqref{knassump}, we have
\begin{equation}
\begin{aligned}\label{betadist}
&\int_0^{k_n-1-l} q_n^{xl}(1-q_n^x)^{l-1}dx=-\frac1{\log q_n}\int_{q_n^{k_n-1-l}}^1 y^{l-1}(1-y)^{l-1}dy\sim\\
&\frac{n^\alpha}c\int_0^1y^{l-1}(1-y)^{l-1}dy=\frac{n^\alpha}c\frac{\Gamma(l)\Gamma(l)}{\Gamma(2l)}=\frac{n^\alpha}c\frac{\big((l-1)!\big)^2}{(2l)!}.
\end{aligned}
\end{equation}
From \eqref{totallowerbound}-\eqref{betadist}, we conclude that
\begin{equation*}\label{finallowermainresult}
P_n^{q_n}(A_{l;k_n}^{(n)})\gtrsim\frac{\big((l-1)!\big)^2}{(2l)!}\frac{c^{l-1}l!}{n^{\alpha(l-1)}},
\end{equation*}
which is the lower bound in  \eqref{mainresultgen} and \eqref{mainresult}.

For the  upper and lower bounds in \eqref{mainresultgen}, the only assumption on $k_n$ was \eqref{knassump}.
It is clear from the proofs that if we fix $\alpha'\in(\alpha,1)$ and let $k_n'=[n^{\alpha'}]$, then the upper and lower bounds in \eqref{mainresultgen} are
uniform over sequences $\{k_n\}_{n=1}^\infty$ satisfying $k_n\ge k_n'$.
From this along with \eqref{duality}, it follows that the upper and lower bounds in \eqref{mainresultgen} are in fact uniform over all sequences $\{k_n\}_{n=1}^\infty$.
\hfill $\square$

\section{Proof  of Theorem \ref{alpha1}}\label{alpha1proof}
 \noindent \it Proof of part (i).\rm\
We follow a slightly more precise version of  the construction used in the upper bound in Theorem \ref{supercluster}, and then reused for the particular case of the Mallows distribution in the proof of Theorem \ref{maintheorem}.
As with the proof of Theorem \ref{maintheorem}, we use the construction from the proof of Theorem \ref{supercluster} in the particular case of the Mallows distribution, with parameter $q_n$; namely, with
$p_j=(1-q_n)q_n^{j-1}$.
Then from \eqref{Xn}, the random variables
$\{X_j\}_{j=2}^\infty$ have truncated geometric distributions.
Although $P_n^{q_n}$ denotes the Mallows distribution with parameter $q_n$, we also use this notation for probabilities of events related to the random variables $\{X_j\}_{j=2}^n$.
It is easy to check that
$P_n^{q_n}(X_m\not\in\{j+1,\cdots,j+l-1\})$ is monotone increasing in $j$. Thus, the  argument leading up to
\eqref{firstupper} in fact gives the following slightly more precise version of \eqref{firstupper}:
\begin{equation}\label{firstupper1}
\begin{aligned}
&P_n^{q_n}(A^{(n)}_{l;k_n}|B_{l;k_n;a})\le \prod_{i=1}^{n-k_n-l+1} P_n^{q_n}(X_{k_n+l-1+i}\not\in\{a+i,\cdots, a+i+l-2\})=\\
&\prod_{i=1}^{n-k_n-l+1}
\big(1-\frac{q_n^{a+i}-q_n^{a+i+l-1}}{1-q_n^{k_n+l+i-1}}\big).
\end{aligned}
\end{equation}
From the assumption on $q_n$,
\begin{equation}\label{sim1}
\prod_{i=1}^{n-k_n-l+1}
\big(1-\frac{q_n^{a+i}-q_n^{a+i+l-1}}{1-q_n^{k_n+l+i-1}}\big)\sim\prod_{i=1}^{n-k_n-l+1}\Big(1-\frac{(l-1)cn^{-1}q_n^{a+i}}{1-q_n^{k_n+l+i-1}}\Big).
\end{equation}
We have
\begin{equation}\label{logsim1}
\begin{aligned}
&\log\prod_{i=1}^{n-k_n-l+1}\Big(1-\frac{(l-1)cn^{-1}q_n^{a+i}}{1-q_n^{k_n+l+i-1}}\Big)=\sum_{i=1}^{n-k_n-l+1}\log\Big(1-\frac{(l-1)cn^{-1}q_n^{a+i}}{1-q_n^{k_n+l+i-1}}\Big)\sim\\
&\int_0^{n-k_n-l+1}\log\Big(1-\frac{(l-1)cn^{-1}q_n^{a+x}}{1-q_n^{k_n+l+x-1}}\Big)dx.
\end{aligned}
\end{equation}
Making the change of variables $y=q_n^x$, using the assumptions on $k_n$ and $q_n$
and defining
\begin{equation}\label{gammacd}
\gamma(c,d)=\log\frac{1-e^{-cd}}{1-e^{-c}}<0,
\end{equation}
in order to simplify notation in the sequel,
we have
\begin{equation}\label{COV1}
\begin{aligned}
&\int_0^{n-k_n-l+1}\log\Big(1-\frac{(l-1)cn^{-1}q_n^{a+x}}{1-q_n^{k_n+l+x-1}}\Big)dx=\\
&-\frac1{\log q_n}\int_{q_n^{n-k_n-l+1}}^1\frac{\log\Big(1-\frac{(l-1)cn^{-1}q_n^ay}{1-q_n^{k_n+l-1}y}\Big)}ydy\le\\
&\frac1{\log q_n}\int_{q_n^{n-k_n-l+1}}^1\frac{(l-1)cn^{-1}q_n^a}{1-q_n^{k_n+l-1}y}dy\sim\\
&-(l-1)q_n^a\int_{e^{-c(1-d)}}^1\frac1{1-e^{-cd}y}dy=\\
&(l-1)q_n^ae^{cd}\log\frac{1-e^{-cd}}{1-e^{-c}}=(l-1)q_n^ae^{cd}\gamma(c,d).
\end{aligned}
\end{equation}
From \eqref{firstupper1}-\eqref{COV1}, we conclude that
\begin{equation}\label{firstupper1final}
P_n^{q_n}(A^{(n)}_{l;k_n}|B_{l;k;a})\lesssim e^{(l-1)q_n^ae^{cd}\gamma(c,d)}.
\end{equation}

Recall that for the particular case of the Mallows distribution with parameter $q_n$, the quantity $N_b$ is given by \eqref{Nb}. Thus, in this particular case, and with $k$ replaced by $k_n$,
\eqref{Ba} becomes
\begin{equation}\label{Baspecific}
P_n^{q_n}(B_{l;k_n;a})=\prod_{j=0}^{l-1}\frac{q_n^a-q_n^{a+j+1}}{1-q_n^{k_n+j}}.
\end{equation}
Using \eqref{firstupper1final} and \eqref{Baspecific}, along with the assumptions on $k_n$ and $q_n$, we have
\begin{equation}\label{upper1final}
\begin{aligned}
&P_n^{q_n}(A^{(n)}_{l;k_n})\lesssim\sum_{a=0}^{k_n-1}q_n^{al}\frac{\prod_{b=1}^l(1-q_n^b)}{\prod_{b=k_n}^{k_n+l-1}(1-q_n^b)}e^{(l-1)q_n^ae^{cd}\gamma(c,d)}\sim\\
&\frac{l!c^l}{n^l(1-e^{-cd})^l}\int_0^{dn}q_n^{xl}e^{(l-1)q_n^xe^{cd}\gamma(c,d)}dx.
\end{aligned}
\end{equation}
Making the change of variables $y=q_n^x$ and using the assumption on $q_n$, we obtain
\begin{equation}\label{COVagain}
\begin{aligned}
&\int_0^{dn}q_n^{xl}e^{(l-1)q_n^xe^{cd}\gamma(c,d)}dx=
-\frac1{\log q_n}\int_{q_n^{dn}}^1 y^{l-1}e^{(l-1)e^{cd}\gamma(c,d)y}dy\sim\\
&\frac nc\int_{e^{-cd}}^1y^{l-1}e^{(l-1)e^{cd}\gamma(c,d)y}dy.
\end{aligned}
\end{equation}
From \eqref{upper1final},  \eqref{COVagain} and \eqref{gammacd},  we arrive at \eqref{alpha1-upper}, which completes the proof of part (i).

\noindent
\it Proof of part (ii).\rm\ We write $q_n=1-\epsilon(n)$, where $0<\epsilon(n)=o(\frac1n)$. We follow the proof of part (i) through the first three lines of \eqref{COV1}, the only change being that
the tern $cn^{-1}$ is replaced by $\epsilon(n)$.
Starting from there, we have
\begin{equation}\label{COVagainn}
\begin{aligned}
&\int_0^{n-k_n-l+1}\log\Big(1-\frac{(l-1)\epsilon(n)q_n^{a+x}}{1-q_n^{k_n+l+x-1}}\Big)dx\le\\
&\frac1{\log q_n}\int_{q_n^{n-k_n-l+1}}^1\frac{(l-1)\epsilon(n)q_n^a}{1-q_n^{k_n+l-1}y}dy=\\
&\frac{(l-1)\epsilon(n)q_n^a}{\log q_n}q_n^{-(k_n+l-1)}\log\big(\frac{1-q_n^n}{1-q_n^{k_n+l-1}}\big).
\end{aligned}
\end{equation}
Since $\epsilon(n)=o(\frac1n)$, we have $1-q_n^n\sim n\epsilon(n)$,
$1-q_n^{k_n+l-1}\sim k_n\epsilon(n)$, $q_n^{-(k_n+l-1)}\sim1$ and
$q_n^a\sim 1$, uniformly over  $a\in\{0,\cdots, k_n-1\}$.
Using this with \eqref{COVagainn}, we have
\begin{equation}\label{COVagainnn}
\begin{aligned}
&\int_0^{n-k_n-l+1}\log\Big(1-\frac{(l-1)\epsilon(n)q_n^{a+x}}{1-q_n^{k_n+l+x-1}}\Big)dx\lesssim (l-1)\log \frac{k_n}n,\\
&\text{uniformly over}\ a\in\{0,\cdots, k_n-1\}.
\end{aligned}
\end{equation}
From \eqref{firstupper1}-\eqref{logsim1} (with $cn^{-1}$ replaced by $\epsilon(n)$) and \eqref{COVagainnn},
we conclude that
\begin{equation}\label{firstupperlittleofinal}
P_n^{q_n}(A^{(n)}_{l;k_n}|B_{l;k;a})\lesssim (\frac{k_n}n)^{l-1}.
\end{equation}
Using \eqref{firstupperlittleofinal} and \eqref{Baspecific}, along with the assumption on  $q_n$, we conclude that
\begin{equation}\label{upperlittleofinal}
\begin{aligned}
&P_n^{q_n}(A^{(n)}_{l;k_n})\lesssim\sum_{a=0}^{k_n-1}q_n^{al}\frac{\prod_{b=1}^l(1-q_n^b)}{\prod_{b=k_n}^{k_n+l-1}(1-q_n^b)}(\frac{k_n}n)^{l-1}\sim\\
&\frac{l!(\epsilon(n))^l}{(k_n\epsilon(n))^l}(\frac{k_n}n)^{l-1}\sum_{a=0}^{k_n-1}q_n^{al}=\frac{l!(\epsilon(n))^l}{(k_n\epsilon(n))^l}(\frac{k_n}n)^{l-1}\frac{1-q_n^{k_nl}}{1-q_n^l}\sim\\
&\frac{l!}{n^{l-1}k_n}\frac{k_nl\epsilon(n)}{l\epsilon(n)}=\frac{l!}{n^{l-1}}.
\end{aligned}
\end{equation}
\hfill $\square$

\section{Proof of Proposition \ref{classif}}\label{classification}
It has already been noted that a  $p$-shifted random permutation with $p_1>0$ and $\sum_{n=1}^\infty np_n<\infty$
satisfies properties (i) and (ii) of the proposition.  From the construction, it is clear that
it also satisfies property (iii),
if the support of the distribution  $p$ is all of $\mathbb{N}$.
Thus, we only need  prove that if a probability distribution $P$ on $S_\infty$
satisfies the three properties stated in the proposition, then it arises as a $p$-shifted permutation for some distribution
$\{p_j\}_{j=1}^\infty$ whose support is all of $\mathbb{N}$ and that satisfies
$\sum_{n=1}^\infty np_n<\infty$.

Let $\Pi$ denote the random permutation under $P$. By property (ii), $\Pi$ is strictly regenerative and the probability
$u_1$ of renewal at the number 1 is positive. (From this it follows that the probability $u_n$ of renewal at the number
$n$ is positive, for all $n$. However, for this proof, we only need the fact that $u_1>0$.) The event
that $n$ is a renewal point, that is, the event $\Pi([n])=[n]$, can be written as
$\cap_{j=1}^\infty \{I_{<n+j}\le j-1\}$.
Thus, we have $u_n=P(\cap_{j=1}^\infty \{I_{<n+j}\le j-1\})>0$. By property (i), this can be rewritten as
\begin{equation}\label{un=}
u_n=\prod_{j=1}^\infty P(I_{<n+j}\le j-1)=\prod_{j=1}^\infty\big(1-P(I_{<n+j}\ge j)\big).
\end{equation}

Recall that the renewal times are labelled
as $\{T_n\}_{n=1}^\infty$. If $n$ is a renewal point, say $T_{k_0}=n$,  then in order that the reduced permutation
$\text{red}(\Pi|_{[T_{k_0+1}]-[T_{k_0}]})$ have the same distribution as $\Pi|_{[T_1]}$,
we need
\begin{equation}\label{I<reduction}
\text{dist}(\{I_{<n+j}\}_{j=1}^\infty|\cap_{j=1}^\infty \{I_{<n+j}\le j-1\})=\text{dist}(\{I_{<j}\}_{j=1}^\infty).
\end{equation}
By property (i), the above reduces
to
\begin{equation}\label{keyreduction}
\text{dist}(I_{<n+j}|I_{<n+j}\le j-1)=\text{dist}(I_{<j}), \ \text{for}\ j=2,3,\cdots\ \text{and}\ n=1,2,\cdots.
\end{equation}
Now the argument leading to \eqref{keyreduction}, for any particular $n$, was arrived at under the assumption that $u_n>0$.
By property (ii), we have $u_1>0$.
Thus, \eqref{keyreduction} holds for $n=1$.
From this it follows that there exist nonnegative $\{p_j\}_{j=1}^\infty$
with $p_1>0$ such that
\begin{equation}\label{p}
P(I_{<j}=i)=\frac{p_{i+1}}{\sum_{k=1}^jp_k},\ i=0,1,\cdots j-1\ \text{and}\ j=2,3,\cdots.
\end{equation}
We now show that $\sum_{j=1}^\infty p_j<\infty$.
Assume to the contrary. Then from \eqref{p} it follows that
$I_{<j}$
converges in probability to $\infty$ as $j\to\infty$.
Thus $\lim_{n\to\infty}P(I_{<n+j}\ge j)=1$, for all $j$ and consequently
$$
\lim_{n\to\infty}\sum_{j=1}^\infty P(I_{<n+j}\ge j)=\infty.
$$
From this and \eqref{un=}, it follows that
$\lim_{n\to\infty}u_n=0$, which contradicts the assumption that the  strictly regenerative random permutation
is positive recurrent.

Since $\sum_{j=1}^\infty p_j<\infty$, without loss of generality we may assume that $\sum_{j=1}^\infty p_j=1$.
From \eqref{p}, we conclude that
$P(I_{<j}=i)=\frac{p_{i+1}}{\sum_{k=1}^jp_k}$, for $i=0,1,\cdots j-1$ and $j=2,3,\cdots$.
From this it follows that the measure $P$ is the $p$-shifted measure with   $p$-distribution given by
$\{p_j\}_{j=1}^\infty$.
In order for property (iii) to hold, it is necessary that $p_j>0$, for all $j$.
\hfill $\square$

\end{document}